\documentclass{article}
\usepackage{amssymb}

\newtheorem{thm}{Theorem}[section] %

\newtheorem{lemma}[thm]{Lemma} %

\newcommand{\eqn}{\begin{eqnarray}}

\newcommand{\eeqn}{\end{eqnarray}}

\usepackage{amsmath}

\begin{document}

\Large{

\centerline{\textbf{The Spectral Geometry of the Mesh Matrices of Graphs}}

\centerline{ Sylvain E. Cappell and Edward Y. Miller}

\centerline{Dedicated to the memory of Jacob Eli Goodman}

\vspace{.3in}

\centerline{\textbf{Abstract}}

}

\vspace{1 in}

The mesh matrix $Mesh(G,T_0)$ of a connected finite graph $G=(V(G),E(G))=(vertices, edges) \ of \ G$ of with respect to a choice of a spanning tree $T_0 \subset G$
is defined and studied. It was introduced by Trent \cite{Trent1},\cite{Trent2}. Its characteristic polynomial  $det(X \cdot Id - Mesh(G,T_0))$ is shown to equal $\Sigma_{j=0}^{N} \ (-1)^j \ ST_{j}(G,T_0)\ (X-1)^{N-j} \  (\star)$ \ where $ST_j(G,T_0)$ is the number of spanning trees of $G$
meeting $E(G-T_0)$ in j edges and $N=|E(G-T_0)|$ . As a consequence, there are  Tutte-type deletion-contraction formulae
for computing this polynomial.   Additionally, $Mesh(G,T_0) -Id$ is of the special form $Y^t \cdot Y$
; so  the eigenvalues of the mesh matrix $Mesh(G,T_0)$ are all real
and are  furthermore be shown to be  $\ge +1$. It is shown that
$Y \cdot Y^t$, called the mesh Laplacian, is a generalization of the standard graph Kirchhoff Laplacian
$\Delta(H)= Deg -Adj$ of a graph $H$.
 For example, $(\star)$ generalizes the all minors matrix tree
theorem for  graphs $H$ and gives a deletion-contraction formula for the characteristic polynomial
of $\Delta(H)$.
This generalization is  explored in some detail. The smallest positive eigenvalue
of the mesh Laplacian, a measure of flux, is estimated, thus  extending the classical
inequality for the Kirchoff Laplacian of graphs.

\newpage

\Large{

\centerline{\textbf{The Spectral Geometry of the Mesh Matrices of Graphs}}

\centerline{ Sylvain E. Cappell and Edward Y. Miller}

\centerline{Dedicated to the memory of Jacob Eli Goodman}

\vspace{.3in}

\centerline{\textbf{Abstract}}

}

\small{

The mesh matrix $Mesh(G,T_0)$ of a connected finite graph $G=(V(G),E(G))=(vertices, edges) \ of \ G$ of with respect to a choice of a spanning tree $T_0 \subset G$
is defined and studied. It was introduced by Trent \cite{Trent1,Trent2}. Its characteristic polynomial  $det(X \cdot Id - Mesh(G,T_0))$ is shown to equal $\Sigma_{j=0}^{N} \ (-1)^j \ ST_{j}(G,T_0)\ (X-1)^{N-j} \  (\star)$ \ where $ST_j(G,T_0)$ is the number of spanning trees of $G$
meeting $E(G-T_0)$ in j edges and $N=|E(G-T_0)|$ . As a consequence, there are  Tutte-type deletion-contraction formulae
for computing this polynomial.   Additionally, $Mesh(G,T_0) -Id$ is of the special form $Y^t \cdot Y$
; so  the eigenvalues of the mesh matrix $Mesh(G,T_0)$ are all real
and are furthermore be shown to be  $\ge +1$. It is shown that
$Y \cdot Y^t$, called the mesh Laplacian, is a generalization of the standard graph Kirchhoff Laplacian
$\Delta(H)= Deg -Adj$ of a graph $H$.
 For example, $(\star)$ generalizes the all minors matrix tree
theorem for  graphs $H$ and gives a deletion-contraction formula for the characteristic polynomial
of $\Delta(H)$.
This generalization is  explored in some detail. The smallest positive eigenvalue
of the mesh Laplacian, a measure of flux, is estimated, thus  extending the classical
inequality for the Kirchoff Laplacian of graphs.

}
\Large{

\section{Introduction}  \label{sect1}
For a  graph $H=(V(H),E(H))$ with vertices $V(H)$ and edges $E(H)$, two distinct vertices, say
$P,Q$, are called adjacent if there is an edge $e$ of $H$ with  $P,Q$ as end points. This is
written as $P \sim Q$. The graph $H$ is called connected if for any two distinct vertices, $P,Q$
there is a simple path from $P$ to $Q$ in $H$, i.e., a sequence of distinct vertices $R_1,R_2, \cdots, R_K$
with $P = R_1,\ Q= R_K$ and $R_k \sim R_{k+1}, \ k=1, \cdots, K-1$. The graph $H$ is called a tree if
$H$ is connected and any
two distinct vertices are connected by a unique simple path. A subgraph $T \subset G$ of a graph $G$
is called a spanning tree of $G$ if the subgraph $T$  is a tree and $V(T)=V(G)$, i.e., $T$ and $G$
contain the same vertices.

\vspace{.1in}  The context of the mesh matrix  studied here is that of a pair $(G,T_0)$ with $G$ a connected finite graph and the chosen subgraph
$T_0 \subset G$  a spanning tree of $G$. It is called $Mesh(G,T_0)$ and is
a real, symmetric, $|E(G-T_0)| \times |E(G-T_0)|$ matrix.

 Before giving the precise definition of the mesh matrix $Mesh(G,T_0)$ of a
 connected graph $G$ with  a choice $ T_0  \subset G$ of spanning tree  of
 $G$ (see section 2) , it will be helpful to recall two relevant items and state two theorems:

\vspace{.1in}
Item 1: An interesting and classically studied invariant of a graph $G$ is the number
of spanning trees  in $G$ \cite{Tutte1}:
$$
ST(G) := \# \ spanning \ trees \ of \ G.
$$
Naturally, $ST(G)=0$ if $G$ is not connected. This invariant satisfies Tutte-type deletion-contraction
formulas.  More generally, if $H \subset G$ is a subgraph and \newline $N= |E(G-H)|$, one may introduce
a mild generalization, the polynomial
$$
ST(G, H)(X) = \Sigma_{j=0}^{N} \ (-1)^j  \ ST_j(G,H) \ X^{N-j}
$$
with $ST_j(G,H)$ the number of spanning trees of $G$ which contain precisely $j$ edges of $E(G-H)$.
Again, see section \ref{sect4}, there are Tutte-type deletion-contraction formulas for
computing these polynomials.

\begin{thm}\label{thmA}  Let $G$ be a connected graph and $T_0 \subset G$ a spanning tree in $G$,
then the characteristic polynomial of the mesh matrix $Mesh(G,T_0)$ is specified by
$$
det( X  \cdot Id - Mesh(G,T_0)) = ST(G, T_0)(X-1).
$$
 See theorem \ref{thm1}
for the complete theorem.
\end{thm}
In view of the above remarks,  this polynomial can be computed directly by deletion-contraction methods.

\vspace{.1in}
Item 2:
It is claimed that $Mesh(G,T_0)$ is of the special form $Id + Y^t \cdot Y$  and
$$
 Y \cdot Y^t, \ where  \ Y \ is \ defined \ via \ (G,T_0),
$$
is a natural generalization of the graph Laplacian $\Delta(H) = Deg - Adj$
of a graph $H$.
This will be  seen as follows:

 Let $H$ be a graph and $W$ be a disjoint  point. Form  $C(V(H))$, the cone on the
vertices $V(H)$ of $H$ from cone point $W$, by attaching an edge $\vec{PW}$ for each vertex $P $ of $H$;
in particular, the vertices of
$C(V(H))$  is the union $V(H) \sqcup \{W\}$. Let the    `` $H$ with cone on vertices added'', be the union
$$
C'(H) = H \sqcup C(V(H)),
$$
i.e., the graph
$C'(H)$ has  vertices $V(C'(H)) = V(H) \sqcup \{W\}$ and edges $E(C'(H)) = E(H) \sqcup
\{ \vec{PW} \ for \ P \ a \ vertex \ of \ H \}$. Note that $C'(H)$ is connected and
$C(V(H))$ is a spanning tree of $C'(H)$.

\begin{thm} \label{thmB}  Let $H$ be a finite graph.
The characteristic polynomials of the Kirchhoff graph Laplacian
$\Delta(H)=Deg - Adj$ of a graph $H$ and that of the mesh matrix of  the pair $(C'(H), C(V(H))$
are related by
$$ \begin{array}{l}
 Mesh(C'(H),C(V(H))) = Id + Y^t \cdot Y \ \ \ with \ \ \
\Delta(H) = Y \cdot Y^t  = \Delta(C'(H),C(V(H))).\\
\end{array} $$
In particular, $\Delta(H)$ has the same non-vanishing eigenvalues counted with multiplicities as
$Y^t \cdot Y = Mesh(C'(H),C(V(H))) -Id$ and as $Y \cdot Y^t = \Delta(C'(H),C(V(H)))$.

See section \ref{sect2}.
\end{thm}

By these means, one sees that the mesh Laplacian defined by  $\Delta(G,T_0)=Y \cdot Y^t$ generalizes to the mesh
setting the well known and studied Kirchhoff graph Laplacian \cite{Kirchhoff1}  of a graph $H$, \  $\Delta(H)$.

 In particular, in view of item \# 1,  the characteristic polynomial
of the graph Laplacian of a graph $H$, namely $det( T \cdot Id - \Delta(H))$,
acquires a Tutte-type deletion-contraction method of computation.

Recently, an elegant paper of  F. Aliniaeifard, V. Wang,and  S. Van Willigenburg, ``Deletion-Contraction for
a Unified Laplacian and Applications'', \cite{Farid1},  has appeared with a deletion-contraction
method for  computing the characteristic polynomial of the graph Laplacian of a graph $H$
among other valuable results.
It proceeds by finding such a formula in terms of an extension $det( T \cdot \epsilon  - \Delta(H))$
where $\epsilon$ is a diagonal matrix of vertex weights. The computation stays in this
restricted  context as opposed to the above which stays in the context defined for
pairs $(G,H)$ with $G$ a graph and $H$ a subgraph. The methods for these cases are analogous,
but different.

\vspace{.1in}

The definition of the mesh matrix $Mesh(G,T_0)$ appears in section \ref{sect2}
along with the proof of theorem \ref{thmB} in a more explicit form in theorem \ref{thm2}.
Section \ref{sect3} contains the proofs of theorems \ref{thmA}, \ref{thm1},  \ref{thmMat}.
In section \ref{sect4} , the Tutte-type deletion-contraction formulas for $ST(G)$ and the polynomial
$ ST(G,H)(X)$ are described.  In section \ref{sect5} ,  by specializing theorem \ref{thm1} to the case $(C'(H),C(V(H)))$ and employing theorem \ref{thmB}, one recovers  the all minors
theorem for graphs; see  theorem \ref{thmminors}.  Section \ref{sect7} proves that $det(Mesh(G,T_0))=ST(G,T_0)$
is the order of the torsion group $C_1(G;\mathbb{Z})/ [ Z_1(G;\mathbb{Z}) \oplus \Pi(B^1(G;\mathbb{Z}))]$
defined there.
Section \ref{sect8} introduces
eigenvalue estimates for the smallest positive eigenvalue of $\Delta(G,T_0)$,
a measure of the flux from $T_0$ to $T_0$ through $G-T_0$. The present authors'
earlier work on mesh matrices on general C-W complexes  appeared in \cite{Cappell1};
section \ref{sect9} below concerns a higher dimensional generalization
of theorem \ref{thm1} with  a correction of a slightly misstated  theorem of  that paper.

\section{Definition of the mesh matrix $Mesh(G,T_0)$:} \label{sect2}
The context of the mesh matrix is that of a pair $(G,T_0)$ with $G$ a connected graph and
$T_0 \subset G$ a spanning tree of $G$.

\textbf{Definition of the mesh matrix : \ $Mesh(G,T_0)$ .}

Now let $T_0$ be a choice of spanning tree of $G$.
Let $e_1, \cdots, e_{|E(G-T_0)|}$ be a listing of edges of $E(G-T_0)$
and $f_1,f_2, \cdots , f_{|E(T_0)|} $ comprise the  edges of $E(T_0)$.

Choose for each edge of $G$ a direction, giving a representation
$\vec{e_j} = \vec{e_j}[P_j, Q_j]$ with $\vec{e}_j$ going from vertex $P_j$ to vertex $Q_j$
and $\vec{f}_k = \vec{f}_k[R_k,S_k]$ going from vertex $R_k$ to vertex $S_k$.

Consider the directed edge $\vec{e}_j \in E(G-T_0)$ with end points $P_j,Q_j$.
Then $P_j, Q_j$ are vertices of the
spanning tree $T_0$; so if $P_j \neq Q_j$ there is a unique directed simple
path from $Q_j$ to $P_j$ in the spanning tree $T_0$, given by distinct vertices  say $Q_j= R_{j,1}, R_{j,2}, \cdots, R_{j,n[j]}= P_j$ with $R_{j,t} \sim R_{j,t+1}$  and the edge
$  [ R_{j,t} \rightarrow R_{j,t+1}]$ chosen in $T_0$.
Let $N= |E(G-T_0)|$.

Now define the mapping
$$ \begin{array}{l}
D : C_1( E(G-T_0); \mathbb{R}) \rightarrow C_1(T_0; \mathbb{R}) \\
\ by  \ \ D( 1 \cdot \vec{e}_j) = 0 \ if  \ 1 \le j \le N\ and \ e_j \ is \ a \ loop \\
\ by \ D( 1 \cdot  \vec{e}_j[P_j,Q_j]) = \Sigma_{t=1}^{n[j]-1} \ \ 1 \cdot \  [ R_{j,t} \rightarrow R_{j,t+1}]
\ \ if \ 1 \le j \le N \ and \ P_j \neq Q_j.
\end{array}
$$
As seen the 1-chain $D( 1 \cdot \vec{e}_j[P_j,Q_k])$ has boundary
$\partial ( D( 1 \cdot \vec{e}_j)) = 1 \cdot P_j - 1 \cdot Q_j$ as the directed
path starts at $Q_j$ and ends at $P_j$.

In this manner one has a uniquely defined 1-chain   for $G$, namely:
$$
Z[j] = 1\cdot \vec{e}_j + D( 1 \cdot  \vec{e}_j)  \ for \ j=1, \cdots,  \ N\ = \ |E(G-T_0)|.
  $$
  This 1-chain is a 1-cycle, it has boundary $0$, and
   $D( 1 \cdot  \vec{e}_j) $ supported in $T_0$.
  If $e_j$ is not a loop, this 1-cycle is a directed simple closed curve, it goes from $P_j$ to $Q_j$ in $E(G-T_0)$, then back through the spanning tree $T_0$ from $Q_j$ to $P_j$.

  Since $Z[j]$ restricted to $E(G-T_0)$ becomes  $ 1\cdot \vec{e}_j $, the linear mapping defined by
  $1 \cdot \vec{e}_j \mapsto Z[j]$ is one to one. It maps $C_1(E(G-T_0); \mathbb{Z})$
  into the integral 1-cycles, $Z_1(G; \mathbb{Z})$. Indeed, this gives an isomorphism:
  $$
  C_1(E(G-T_0); \mathbb{Z}) \stackrel{\cong}{\rightarrow} Z_1(G; \mathbb{Z}).
  $$
In particular, these 1-cycles $\{Z[j]\}$ form an integral basis for the integral 1-cycles of $G$.

  [Proof: For any  1-cycle $z= \Sigma_{j=1}^{N} \ b_j \cdot \vec{e}_j+
   \Sigma_{k=1}^{|E(T_0)|} \ c_k \cdot \vec{f}_k$ of $G$,  the 1-chain
  $
  z -\Sigma_{j=1}^N  \ b_j \cdot Z[j]
  $
  is supported in $T_0$ and is a 1-cycle, so equals the zero chain.] The mesh matrix of
  Trent \cite{Trent1,Trent2} may be  defined using this integral basis,  $\{Z[j]\}$, of 1-cycles of $G$.

  Having chosen a spanning tree $T_0$ of the graph $G$, and having made these directed edge choices, one defines the
  mesh matrix for the pair $(G,T_0) $ via the intersection pairing on $C_1(G; \mathbb{R})$:
  $$
  Mesh(G,T_0) = \{\  < Z[j_1], \ Z[j_2]>\ | \ 1 \le j_1,j_2 \le N = |E(G-T_0)|  \}.
$$

The basic theorem about the mesh matrix $Mesh(G,T_0)$ is the theorem of Trent  evaluating
its determinant. This is studied in more detail in \cite{Cappell1} along with higher dimensional
generalizations.

\begin{thm} [Trent \cite{Trent1,Trent2} ]\label{thmTrent1} For a connected graph $G$, the determinant \newline $det(Mesh(G,T_0)) $ equals the number
of spanning trees  of $G$,  or as in item \# 1 equals $ST(G)$. \end{thm}

 More explicitly,  expand the 1-cycle  $Z[k]$ as a column $|E(G)| \times 1$ vector, say $z[k]$,
in terms of the basis $\{1 \cdot \vec{e}_j\}$ and then $\{1 \cdot \vec{f}_k\}$ of $C_1(G;\mathbb{Z})$
and form the $ |E(G)| \times N$ matrix of column vectors $Y[+] := ( z[1],\cdots, z[N])$, then
$
Mesh(G,T_0) = Y[+]^t \cdot Y[+],
$
As seen $Y[+]$ has the special form
$
Y[+] = \left( \begin{array}{c} \ Id \\ Y \end{array} \right)
$
where $Id$ denotes the $N \times N$ identity matrix and $Y$ is the $|E(T_0)| \times |E(G-T_0)|$ matrix representing  the linear mapping $D$.

As a consequence,
$$
\begin{array}{l}
Mesh(G,T_0) = Y[+]^t \cdot Y[+] = Id \ \ + \ \ (Y^t \cdot Y) \  and \ \\
charactistic \ polynomial \ of \ Mesh(G,T_0) = det(X \ Id - Y[+]^t \cdot Y[+])
  \\ = det((X-1) \ Id - Y^t \cdot Y)\\ = characteristic \ polynomial \ of \ (Y^t \cdot Y) \ evaluated \ at  \ (X-1).
\end{array}  $$

In particular, the eigenvalues of the reduced mesh matrix $Mesh^\#(G,T_0) \stackrel{def.}{:=} Y^t \cdot Y$ are real and non-negative, while
those of the mesh matrix $Mesh(G,T_0)$ are obtained by adding $+1$ to those of the reduced mesh matrix.

The following theorem specifies the coefficients of the characteristic polynomial of the mesh matrix and the reduced mesh matrix of $(G,T_0)$. It is a restatement of theorem \ref{thmA}.

 \begin{thm} \label{thm1} Let $G$ be a connected graph and $T_0$  a spanning tree of $G$. For $j=1, \cdots, N= |E(G-T_0)|$, denote by $ST_j(G,F)$
 the number of spanning trees $T$ of $G$ such that $|E(T) \cap (E(G-T_0)) |= j$. Then
 $$  \begin{array}{l}
 det((U+1)  Id - Mesh(G,T_0))= det( U \ Id - Mesh^\#(G,T_0)) \\ = U^N + \Sigma_{j=1}^N \ (-1)^j  ST_j(G,T_0)\  U^{N-j} = ST(G,T_0)(U).  \end{array}
 $$
Setting $U=-1$ gives $det(Mesh(G,T_0)) = \Sigma_{j=0}^N \ ST_j(G,T_0) = ST(G)$ in agreement
with Trent's theorem \ref{thmTrent1}. Here $Y^t \cdot Y$ represents the reduced mesh matrix $Mesh^\#(G,T_0)$ and $Mesh(G,T_0)= Id + Y^t \cdot Y$.

 \end{thm}

In terms of the previously defined polynomials for $H$ a subgraph of $G$ listed in item \# 1, section \ref{sect1}
$
ST(G,H)(U),
$
this theorem evaluates the characteristic polynomials of the  reduced mesh matrix $Y^t \cdot Y$
as equal to  $ST(G,T_0)(U)$. In this larger context of pairs of graphs, $(G,H)$ one has, see section \ref{sect4},
deletion-contraction formulas which  may be used to compute this  invariant.

 \vspace{.3in}
 It is natural to define the $|E(T_0)| \times |E(T_0)|$ matrix
 $$
 Y \cdot Y^t  \ as \ the  \  Mesh \ Laplacian\ of \ (G,T_0) \ = \Delta(G,T_0)= D \cdot D^\star.
 $$
 This terminology is consistent with the standard notations in view of the following theorem.
 [See \cite{Chung1,Chung2} for a discussion of $\Delta(G)= Deg-Adj$ and its normalized version.]

 As in item \# 2, section \ref{sect1}, for a graph $H$, let $C'(H)= H \cup C(V(H))$ be $H$ with the  cone on vertices of    $H$ added,

\begin{thm} \label{thm2} For a  graph $H$, consider  $C'(H)$
and its subcone  $C(V(H))$ with cone point $W$. Then $C'(H)$ is connected and $C(V(H))$ is a spanning
tree for $C'(H)$.

Following the above construction, for each edge of $H$ once directed, say \newline
$\vec{e}_j= \vec{e}_j [P_j \rightarrow Q_j], j=1, \cdots, |E(H)|$   associate the 1-cycle, $Z[j] \newline =1\cdot \vec{e}_j
+ 1\cdot (Q_j \rightarrow W) - 1 \cdot (P_j \rightarrow  W)$ of $C'(H)$, then the associated
 matrix $Y$ defined by
 $$
 \begin{array}{l}
 1 \cdot \vec{f}_j \mapsto D(1 \cdot \vec{f}_j) = 1\cdot (Q_j \rightarrow W) - 1 \cdot (P_j \rightarrow W)
 \in C_1(C(V(H));Z) \\
 under \ the \ identification \ C_1(C(V(H));Z) \cong C_0(V(H);Z)\  given \ by \\
  \ 1 \cdot (P \rightarrow W)
 \mapsto 1 \cdot P
 \end{array}
 $$
 becomes [under this identification] precisely the standard  boundary  mapping \newline  $\partial : C_1(H;Z) \rightarrow C_0(H;Z)$,
 $\partial(1\cdot \vec{f}_j ) = 1 \cdot Q_j - 1\cdot P_j$ of the graph $H$.

 Consequently, the mesh Laplacian  of $(G,T_0)= \Delta(G,T_0) $ represented by  $ Y \cdot Y^t $ for this special case  $(C'(H),C(V(H))$
 becomes exactly the graph Laplacian matrix $\Delta(H) =\partial \cdot \partial^t= Deg - Adj$ of $H$.

Also $Mesh(C'(H),C(V(H)) = Id + Y^t \cdot Y = Id + Mesh^\#(C'(H),C(V(H))$.
 \end{thm}

 The significance of this theorem is that the mesh Laplacian matrix $ Y \cdot Y^t$
 representing $ D \cdot D^\star$
 may be viewed as a natural generalization of the Kirchhoff graph Laplacian $\Delta(H)$ to the context of
 mesh matrices for pairs $(G,T_0)$.

  Since eigenvalues of $Y \cdot Y^t$ and $Y^t \cdot Y$
 are non-negative, and their non-vanishing eigenvalues are equal [counted with multiplicities]
 the characteristic polynomials of $Y \cdot Y^t$ and $= Y^t \cdot Y $
 are related by
 $$
  U^{|E(T_0)|} \cdot det(U \ Id - Y^t \cdot Y) =  \\
  U^{|E(G-T_0)|} \cdot det(U \ Id - Y \cdot Y^t).
  $$
  Applied to the case of  $C'(H)$ theorem \ref{thmB} of item \# 2 follows.

 Additionally, the mesh Laplacian $\Delta(G,T_0) = Y \cdot Y^t$ has a nice direct description:

 \begin{thm} \label{thmMat} Let the edges of $E(T_0)$ be directed and labeled as $\vec{f}_k, k=1, \cdots, |E(T_0)|$.
 Then the mesh Laplacian $\Delta(G,T_0) = Y \cdot Y^t$  is specified by
 $$ \begin{array}{l}
 ( Y \cdot Y^t)_{\vec{f}_k, \vec{f_k}} = \# \{ e \in E(G-T_0) \ with \ f_k \ in \ the \ support \ of \ D(1 \cdot \vec{e}) \},  \\
 for \ k \neq l \ \\
 ( Y \cdot Y^t)_{\vec{f}_k, \vec{f_l}}  = Sign(k,l) \cdot \# \{ e \in E(G-T_0) \ with \ f_k \
 and \  f_l \  in \ the \ support \ of \ D(1 \cdot \vec{e}) \}
 \end{array}
 $$
 Here $Sign(k,l)= \pm 1$ for $k \neq l$ is defined by taking the unique directed path, $\gamma$,  in $T_0$
 from the edge $f_k$ to the edge $f_l$.  Let $Sign(k,l) = < 1 \cdot \vec{f}_k, \vec{\gamma}> \cdot <1 \cdot  \vec{f}_l , \vec{\gamma}>$.
 \end{thm}

\section{Proof of theorems \ref{thmA}, \ref{thm1},  \ref{thmMat} : } \label{sect3}
                  Now the characteristic polynomial of $Mesh(G,T_0)$ represented by $Y[+]^t \cdot Y[+]$ is  easily
                  determined by the following method: Let
                  $$
                  det( X \ Id - Y[+]^t \cdot Y[+]) = X^N + \Sigma_{j=1}^N (-1)^j \ b_j \ X^{N-j}.
                  $$

                  Firstly, by   expanding to get the coefficient of $X^{N-j}$,
                  $b_j$ equals the sum over $j \times j$ \textbf{diagonal minors}  of $Y[+]^t \cdot Y[+]$ of
                  the determinant  $det(minor)$.
                  This sum may be organized as summing over a choice of $j$ indices
                  say $1 \le k_1 < k_2 < \cdots < k_j \le N$ [specifying the rows and corresponding
                  columns  of the diagonal block chosen].

More explicitly,   let $Y[+][k_1,\cdots, k_j]$ denote the $|E(G)|  \times j$ matrix
from $Y[+]$ obtained
                   by deleting all but these $j$ specified columns. Then the desired block  determinant
                   equals
                   $$
                   det( \ Y[+][k_1,\cdots, k_j]^t \cdot Y[+][k_1,\cdots, k_j] \ ).
                   $$

                   Now compare this last  to the determinant of the mesh matrix for
                   the subgraph \newline  $G[k_1,k_2, \cdots, k_j]$ obtained from the spanning
                   tree $T_0$ by adding in the $j$ edges $\vec{e_{k_a}}, \ a=1, \cdots, j$
                   of $E(G-T_0)$ to $T_0$. Each of the 1-cycles $Z[k_a]$ lie in $G[k_1,k_2, \cdots, k_j]$. As seen, $ Y[+][k_1,\cdots, k_j]^t \cdot Y[+][k_1,\cdots, k_j]
                   = Mesh(G[k_1,k_2, \cdots, k_j];T_0)$. Consequently, by
                   the results of Trent:
                   $$   \# spanning \ trees  \ of \ G[k_1, \cdots, k_j]
                    = det( \ Y[+][k_1, \cdots, k_j]^t \cdot  Y[+][k_1, \cdots, k_j] \ ).
                    $$

                    So in total,  the coefficient $b_j$ is expressed  in terms
                    of spanning trees. More explicitly :

                    \begin{lemma} \label{lemmaA} Let $G$ be a connected graph and $T_0$ be a spanning tree
                    of $G$ and $N= |E(G-T_0)|$. Then the characteristic
                    polynomial of the mesh matrix $Mesh(G,T_0)$ is given by
                    $$
                  det( X \ Id - (Y[+]^t \cdot Y[+]) =X^N+ \Sigma_{j=0}^N (-1)^j \ b_j \ X^{N-j}
                  $$
                   with $b_j$ equal to the sum over $j$ distinct indices
                   $1 \le k_1 < k_2< \cdots k_j \le N$ of the number of spanning
                   trees  in the subgraph $ T_0 \cup ( e_{k_1} \cup e_{k_2} \cup \cdots \cup e_{k_j}) $.
                   \end{lemma}

                   \vspace{.3in}

   Now consider this sum $b_j$ in detail with reference to the integers
   $$
   ST_j(M,T_0) = \# \{ \ spanning \ trees \ T \ of \ G \ with \ |E(T) \cap E(G-T_0)| = j \ \}.
   $$

   Each spanning tree $T$ counted in $ST_j(G,T_0)$ appears exactly once in the counting
   of $b_j$ as it fills in all the j added edges of $E(G-T_0)$.

   Each spanning tree $T$ of $ST_{j-1}(G,T_0)$ fills in  $j-1$ edges of $E(G-T_0)$
   so to get $j$ one  needs to chose 1 from the remaining $N-(j-1)$ edges. Thus
   each spanning tree  $T$ of $ST_{j-1}(G,T_0)$  appears exactly $ {{N-(j-1)}\choose{1}}
   = {{N-(j-1)}\choose{N-j}}
   $ times in $b_j$.

Each spanning tree  $T$ of $ST_{j-2}(G,T_0)$ fills in  $j-2$ edges of $E(G-T_0)$
   so to get $j$ one  needs to chose 2 from the remaining $N-(j-2)$ edges. Thus
   each spanning tree $T$ of $ST_{j-2}(G,T_0)$  appears exactly $ {{N-(j-2)}\choose{2}}
   = {{N-(j-2)}\choose{N-j}}
$ times in $b_j$.

Combining these observations, proves the formula for $ N \ge j \ge 1$:
$$
b_j  = \Sigma_{t=0}^j \ {{N-(j-t)}\choose{N-j}} \ ST_{j-t}(G,T_0).
$$
It also holds for $j=0$ since both sides equal +1. Inserting [and letting $k=j-t$] this gives:
$$\begin{array}{l}
det( X \ Id - Mesh(G,T_0)) = \Sigma_{j=0}^{N} \ (-1)^j \  b_j \ X^{N-j} \\
= \Sigma_{j=0}^{N} \ (-1)^j [ \Sigma_{t=0}^j \ {{N-(j-t)}\choose{N-j}} \ ST_{j-t}(G,T_0)]  \ X^{N-j}\\
= \Sigma_{k=0}^{N} ( \ \Sigma_{j=k}^N \ (-1)^{j} \ {{N-k}\choose{N-j}} X^{N-j}\ ) \ ST_k(G,T_0  ) \\
=  \Sigma_{k=0}^N  \ (-1)^k \ ( \Sigma_{j=k}^{N} (-1)^{j+k} \ {{N-k}\choose{N-j}} X^{N-j}) \ ST_k(G,T_0) \\
=  \Sigma_{k=0}^N  \  (-1)^k \ ( X-1)^{N-k} \ ST_k(G,T_0) = det((X-1) \ Id - Mesh^\#(G,T_0)) . \end{array}
$$
Thus $$ \begin{array}{l} det( (U+1) \ Id - Mesh(G,T_0))=det( U \ Id - Mesh^\#(G,T_0))\\ = U^N + \Sigma_{j=1}^N \ (-1)^j \ ST_j(G,T_0) \ U^{N-j} \end{array} $$
as claimed.

\vspace{.1in}
Proof of theorem \ref{thmMat}:

For directed edge $e_j$  with end points directed by $P_j \rightarrow Q_j$
let $C[Q_j,P_j]$ be the unique simple path from $Q_j$ to $P_j$ in $T_0$.
By definition $D( 1 \cdot \vec{e}_j[P_j,Q_j])$ is then the sum of these
directed edges of $C[P_j,Q_j]$ regarded as a 1-chain.

Tautologically, the  mapping $D$ and adjoint $D^\star$ have the properties:
$$  \begin{array}{l}
If \ P_j \neq Q_j, \ D( 1 \cdot \vec{e}_j[P_j,Q_j]) = \Sigma_{f_k \ in \ support \ of \ C[Q_j,P_j]}
\ \  < 1 \cdot \vec{f}_k, D( 1 \cdot \vec{e}_j[P_j,Q_j]) > \cdot \vec{f}_k; \\
D( 1 \cdot \vec{e}_j[P_j, Q_j]) = 0 \ if \ e_j \ is \ a \ loop \ i.e., P_j=Q_j.  \\

Similarly,\\
if \ f_k \ not \ a \ loop,\ \ D^\star( 1 \cdot \vec{f}_k) = \Sigma_{f_k \ in \ support \ of \ C[Q_j,P_j]} \
\  < 1 \cdot \vec{f}_k, D( 1 \cdot \vec{e}_j[P_j,Q_j])> \cdot \vec{e}_j ; \\
D^\star( 1 \cdot \vec{f}_k) = 0 \ if \ f_k \ is \ a \ loop . \\\end{array} $$

Hence, by $< f, \Delta(G,T_0)(f)>= < D^\star(f), D^\star(f)>$, the mesh Laplacian  matrix is defined by
$$ \begin{array}{l}
\Delta(G,T_0)_{f_k,f_l} = \Sigma_{e_j \ with \ f_k\ and \ f_k \ in \ the \ support \ of \
C[Q_j,P_j]} \\  \hspace{1in}  < 1 \cdot \vec{f}_k, D( 1 \cdot \vec{e}_j[P_j,Q_j])> \cdot < 1 \cdot \vec{f}_l, D( 1 \cdot \vec{e}_j[P_j,Q_j])>.
\end{array}
$$
For $k=l$,  this give the desired result. For $k \neq l$,  one may
change the ordering of $Q_j,P_j$ if necessary so that the unique  path $C[Q_j,P_j]$
running from $Q_j$ to $P_j$ in $T_0$ containing the edges $f_k,f_l$ gives the unique directed
path, say $\vec{\gamma}$,  in $T_0$ from $f_k$ to $f_l$. In particular, for any $e_j$ with
$f_k,f_l$ in $C[Q_j,P_j]$  one has the equality independent of $j$
$$  \begin{array}{l}
 < 1 \cdot \vec{f}_k, D( 1 \cdot \vec{e}_j[P_j,Q_j])> \cdot < 1 \cdot \vec{f}_l, D( 1 \cdot \vec{e}_j[P_j,Q_j])>
 = < 1 \cdot \vec{f}_k,\vec{ \gamma}> \cdot < 1 \cdot \vec{f}_l, \vec{\gamma}> , \ \\
 \end{array} $$
 as desired.

\section{Properties of $ST(G,H)(X)$; Computation by deletion contraction methods:} \label{sect4}
 It is well known that the number of spanning trees, $SF(G)$, obeys nice formulas
 for deletion  and contraction \cite{Tutte1}.

 For a graph $G$ and a edge $e$, let $G\backslash e$ be the graph obtained by
 deleting the open edge $e$ preserving the end points, i.e, $V(G\backslash e) = V(G)$
  and $E(G\backslash e) = E(G) - \{e\}$.
 Let $G / e$ denote the
 graph obtained by contracting the end points of $e$ along with  all of the edge $e$
 to a single  point.

 As $ST(G)$ counts  spanning trees, $ST(G) = 0
 \ if \ G \ is\  disconnected. $

 As no tree contains a loop,
 $
 ST(G) = ST(G \backslash e) \ if \ e \ is \ a \ loop.
 $

 If $e$ is an edge with distinct end points, then the spanning trees of $G$ counted in $SF(G)$ either
 do not contain $e$ and are faithfully counted in $ST(G \backslash e)$ or contain the edge $e$.
 In this second case, the spanning trees of $G/e$ are exactly the result of collapsing the
 edge $e$ in a spanning tree of $G$ containing $e$. In this manner, one gets the well
 known relations:
 $$
 ST(G) = ST(G\backslash e) + ST( G / e) \ if \ e \ an \ edge \ of \ G \ is \ not \ a \ loop.
 $$

 If $e$ is an isthmus, that is all trees of $G$ contain the edge $e$, then $G \backslash e$ is disconnected, so
 $$
 ST(G) = ST(G\backslash e) + ST( G / e) = ST(G/e) \ if \ e\ is \ an  \ isthmus.
 $$

 By using these relations in the context of a graph $G$ and a subgraph $H$ with
 $$
 ST(G,H)(X) = \Sigma_{j=0}^{N} \ (-1)^j  \ ST_j(G,H) \ X^{N-j} \ and \ N =|E(G-H)|
 $$
 where $ST_j(G,H)$ is the number of spanning trees $T$ of $G$ with $|E(T) \cap (E(G)-E(H))| = j$,
 one obtains by counting:
 $$
 \begin{array}{l}
 ST(G,H)(X) = ST(G ,H \backslash e)(X) \ if \ e \in H \ is \ a \ loop; \\
  ST(G,H)(X) = ST(G ,H \backslash e)(X)+ ST(G ,H / e)(X) \ if \ e \in H \ is \ not\ a \ loop; \\
ST(G,H)(X) = ST(G \backslash e ,H )(X) \ if \ e \in G \ is \ a \ loop; \\
ST(G,H)(X) = ST(G \backslash e,H )(X)- X \cdot ST(G / e ,H )(X) \\ \hspace{.5in} if \ e \in G\ is \ not \ a \ loop.
\end{array}
$$

These deletion-contraction formulas give an inductive method for computing
the polynomials $SF(G,H)$.

By theorem \ref{thm1}, $det(X \cdot Id - Mesh^\#(G,T_0)) = ST(G,T_0)(X)$
so these formulas give a deletion-contraction method for computing the
characteristic polynomial of the reduced mesh matrix $Mesh^\#(G,T_0)$
and the mesh matrix $Mesh(G,T_0)= Mesh^\#(G,T_0) + Id$.

As a special case, by theorem \ref{thm2}, the characteristic polynomial of the graph Laplacian
$\Delta(H)$ of a graph $H$ equals that of $\Delta(C'(H),C(V(H)))$ so these deletion-contraction
formulas give a method that also applies to the graph Laplacian $\Delta(H)$.

Recall from \ item \# 2, section \ref{sect1}, the paper \cite{Farid1} gives
an alternative method to get deletion-contraction formulas for computation in these graph Laplacian cases.

\section{Recapturing the All Minors Matrix Tree Theorem} \ \label{sect5}

The all minors matrix tree theorem, see \cite{Chaiken1,Chaiken2,Maurer1}, asserts:
\begin{thm}  \label{thmminor1} Let $H$ be a finite   graph. Let a spanning forest $F$  in $H$
consist of the union of $K$ disjoint trees, say $F= \sqcup_{j=1}^K \ T[j],  \subset H$ which spans,
i.e., $V(F)= V(H)$.
Let a spanning  forest $F$ of $K$ components have multiplicity $mult(F) = \prod_{j=1}^K
\ |V(T[j])|$. Then the characteristic polynomial of the graph Laplacian $\Delta(H)
= \partial \cdot \partial^\star = Deg - Adj $ is given by
$$
\begin{array}{l}
det( T \cdot Id - \Delta(H))=T^{|V(H)|}+  \Sigma_{j=1}^{|V(H)|- 1}    \ (-1)^{j} \ b_j(H) \              T^{|V(H)| - j} \\
with \ b_j(H) = \Sigma_{F \ a \ spanning \ forest \ of \ H\ with \ |V(H)|-j \ components} \
 \ mult(F).
\end{array}
$$

 Note: Let a spanning rooted forest consist of a spanning forest $F$ and a choice of a vertex, say  $P_j \in V(T[j])$, called roots
for each $j$.   Then $$b_j(H) = \#  \ spanning \ rooted\  forests \ in \ H  \ with \ |V(H)|-j \ components$$ .
\end{thm}

Now by contrast consider the spanning trees $T$ of $C'(H)= H \cup C(V(H))$,
  obtained from
the graph $H$  with the cone on the vertices added. Necessarily, the spanning tree $T$ must contains the cone point,  say $W$, and
all vertices of $H$. Also, the intersection $F = T \cap H$ is a spanning  forest in $H$,
say $F = \sqcup_{j=1}^K \ T[j]$ with $K$ components and each $T[j]$ a  tree in $H$. To complete the tree $T$ in $C'(H)$  each component of $F = \sqcup_{j=1}^K \ T[j]$ must be joined by exactly one edge of the form
$\vec{P_jW}$ with $ P_j \in V(T[j])$. So one has the equality

$$ \begin{array}{l}
\# \{ spanning \ trees \ T \  in \ C'(H) \ with \ |E(T) \cap C(V(H))| = K\} \\
\# \{  spanning \ rooted \ forests \ F \ in \ H \ with \ K \ components \} \\
=\Sigma_{spanning \ forests \ F \  in \ H \ with \ K \ components} \ \  mult(F) .
\end{array}
$$

By these means the all minors matrix tree theorem can be recast in the form
\begin{thm} [All minors matrix tree theorem] \label{thmminors}:
$$
\begin{array}{l}
det( T \cdot Id - \Delta(H))=T^{|V(H)|}+  \Sigma_{j=1}^{|V(H)|- 1}    \ (-1)^{j} \ b_j(H) \             T^{|V(H)| - j} \\
with \ b_j(H) = \# \{ spanning \ trees \ T \  in \ C'(H) \ with \ |E(T) \cap C(V(H))| =|V(H)| -j \} \\ = \# \{ spanning \ trees \ T \  in \ C'(H) \ with \ |E(T) \cap E(H)| =j \} \\
= ST_j(C'(H),C(V(H)))
\end{array}
$$
The last equalities use the fact that a spanning tree $T$ of $C'(H)$ will have $|V(C'(H))| - 1
=|V(H)|$ edges.
\end{thm}

Hence, one sees that the all minors matrix tree theorem
is the special case with $(G,T_0) = (C'(H), C(V(H)))$ of theorem \ref{thm1}
in view of theorem \ref{thm2}.

\section{ \# spanning trees as the order of a torsion group associated to $G$:} \label{sect7}

Let $G$ be a connected graph with choices of directed edges as above.  Let $Z_1(G;Z)$ be the integral 1-cycles, i.e.,
the kernel of the boundary mapping $ \partial : C_1(G;\mathbb{Z}) \rightarrow C_0(G;\mathbb{Z})$;
let $B^1(G;Z)$ be the 1-coboundaries i.e., the image of the dual coboundary mapping
$\partial^\star : C^0(G;\mathbb{Z})= Hom(C_0(G;\mathbb{Z}),\mathbb{Z}) \stackrel{Hom(.,\partial)}{\rightarrow} C^1(G;\mathbb{Z})=Hom(C_1(G;\mathbb{Z}), \mathbb{Z})$; let
$$
\begin{array}{l}
\Pi : C^1(G;\mathbb{Z}) \rightarrow C_1(G;\mathbb{Z}) \ be \ the \ isomorphism\ defined \ by \\
  {[ F: E(G) \rightarrow \mathbb{Z}]} \mapsto [ \Sigma_{\sigma \in E(G)} \ F(\sigma) \cdot \vec{\sigma}  ] \\
  \ with \ chosen \ direction  \ \vec{\sigma} \  on \ an  \ edge \ \sigma \ of \ G.
  \end{array} $$

\begin{thm} \label{thmtorsion} The inclusion of integral chains
$
[Z_1(G;\mathbb{Z}) \oplus \Pi( B^1(G;\mathbb{Z}))] \subset C_1(G;\mathbb{Z})
$
is a rational equivalence; so the quotient
$
U = C_1(G;\mathbb{Z})\ / \ [Z_1(G;\mathbb{Z}) \oplus \Pi( B^1(G;\mathbb{Z}))]
$
is a finite abelian group. The order of this group is
$$
| \ U \ | = | \ C_1(G;\mathbb{Z})\ / \ [Z_1(G;\mathbb{Z}) \oplus \Pi( B^1(G;\mathbb{Z}))] \ | = \# \ spanning \ trees \ of \ G
=ST(G)
$$
\end{thm}

First,  observe that the orthogonal complement of $Z_1(G;\mathbb{R})$  is $\Pi( B^1(G;\mathbb{R}))$,
so one  gets an isomorphism of real vector spaces
$
Z_1(G;\mathbb{R}) \oplus \Pi(B^1(X;\mathbb{R})) = C_1(G;\mathbb{R}).
$

Proof:
If $<.,.>$ is the standard inner product on $C_1(X;\mathbb{R})$, then one has  for \newline $a \in Z_1(G,\mathbb{R}), b \in C_0(G;\mathbb{R})$,
$
< a, \Pi( \delta(b))> = \delta(b)(a) = b ( \partial a).
$
Hence,  $\partial a  =0 $ if and only if $a$ is orthogonal to $\Pi(B^1(G;\mathbb{R})))$.

By this result,  the associated embedding of lattices has  finite quotient,
$$
U = C_1(G;\mathbb{Z})\ / \ [Z_1(G;\mathbb{Z}) \oplus \Pi( B^1(G;\mathbb{Z}))].
$$

Next, let $T_0$ be a choice of  spanning tree for $G$. Following the notation of section \ref{sect2}, let the directed edges of $G-T_0$ be denoted by $\vec{e}_j,
j=1,\cdots, |E(G-T_0)|$. Let the directed edges of $T_0$ be denoted by $\vec{f}_k, k=1, \cdots ,|E(T_0)|$.

As above, consider the integer 1-cycles $Z[j], j=1, \cdots, |E(G-T_0)|$.  Here the 1-cycle
$$
Z[j] = 1 \cdot \vec{e}_j + D(\vec{e}_j), \
 j=1, \cdots, |E(G-T_0|
 $$
 has $D(\vec{e}_j)$ supported in the tree $T_0$. As above, define the $ |E(T_0)| \times |E(G-T_0)| $ matrix
 $Y= \{ A[k,j]\}$  by
 $$
 Z[j] = 1 \cdot \vec{e}_j + \Sigma_{k=1}^{|E(T_0)|} \ Y[k,j] \ \vec{f}_k, \ \ \ j=1,\cdots, |E(G-T_0)| .
 $$
Here the coefficients $Y[k,j] $ are either $0$ or $\pm 1$. As proved above, these $Z[j]$'s form
an integral basis for $Z_1(G;\mathbb{Z})$. That is, part a) below is known.

\begin{lemma} \label{lemma1} a)  The 1-chains
$$
Z[\vec{e}_j]= 1 \cdot \vec{e}_j  +D( \vec{e}_j)  =  \vec{e}_j + \Sigma_{k=1}^{|E(T_0)|} \ A[k,j] \ \vec{f}_k , \  \ \ j=1, \cdots, |E(G-T_0)|
$$
form an integral basis for $Z_1(G;\mathbb{Z})$.

b) The 1-chains
$$
b[\vec{f}_k] = 1 \cdot \vec{f}_k  - \Sigma_{j=1}^{|E(G-T_0)|}  \ Y[j,k]\ \vec{e}_j,  \ \
k =1, \cdots, |E(T_0)|
$$
form an integral basis for $\Pi(B^1(G;\mathbb{Z}))$.
\end{lemma}

Using the direct sum  decomposition
$
C_1(G;\mathbb{Z}) = C_1(G-T_0;\mathbb{Z})  \oplus C_1(T_0;\mathbb{Z}),
$
write elements in column form:
$
C_1(G;\mathbb{Z}) = \left( \begin{array}{cc} C_1(G-T_0;\mathbb{Z}) \\ C_1(T_0;\mathbb{Z}) \end{array} \right) .
 $

As noted,  the $\vec{e}_j \mapsto Z[\vec{e}_j]$ giving an integral basis of $Z_1(G;Z)$
is represented by the matrix
$
  \left( \begin{array}{cc} Id \\ Y \end{array} \right)
  $
  with $Y$ an integer valued matrix above. Granting lemma \ref{lemma1} part b),
 the mapping $\vec{f}_k \mapsto b[\vec{f}_k]$ giving an integral basis of $\Pi(B^1(G;\mathbb{Z}))$
  is represented by the matrix
  $
  \left( \begin{array}{cc} -Y^t \\ Id \end{array} \right).
$

Thus,
$
|U|  = det \left( \begin{array}{cc} Id & - Y^t \\ Y & Id \end{array} \right)
= det \left( \begin{array}{cc} Id + Y^t \ Y& - Y^t +Y^t\\ Y  & Id \end{array} \right)
= det(Id+ Y^t \ Y) =det(Mesh(G,T_0)=  ST(G)$  by Trent's theorem \ref{thmTrent1}.

\vspace{.1in}
Proof of lemma \ref{lemma1} part b):

The first part about the 1-cycles $Z[j]$ has already been proved in section \ref{sect2}.

Pick a directed edge in $T_0$, say $\vec{f}_k  = \vec{f}_k[R_k \rightarrow S_k]$. Recall that
the spanning tree $T_0$ is connected.

Decompose the set of vertices of $T_0$ into two disjoint non-empty sets
$
V(T_0)= A[k] \sqcup B[k]
$
where $A[k]$ consists of the vertices of $T_0$ which can be connected to the initial vertex $R_k$
by a simple path not including $f_k$ and $B[k]$ consists of the vertices of $T_0$ which can be connected to the final vertex $S_k$
by a simple path not including $f_k$. For convenience, chose the directions of the edges $\vec{e}_j$ [
which by assumption range over $E(G-T_0)$]
so that if $\vec{e}_j$ has one end point in $A[k]$ and the other in $B[k]$,
then $\vec{e}_j$ is directed from the point in $A[k]$ to the point in $B[k]$.

Now recall that for a vertex $P$ of $G$, the coboundary $\delta(1 \cdot P)$ has the property that
$\Pi(\delta(X))$ is the 1-chain consisting of all edges with distinct end points containing
$P$ as a end point and directed [inwards] towards $P$.

One easily sees  that $\Pi(\delta( \Sigma_{b \in B[k]} \ 1 \cdot b ))$
consists of a sum of two parts. Firstly $1 \cdot \vec{f}_k$; and secondly a sum
over all edges $\vec{e}_j$ in $G-T_0$  with distinct end points and having
one end in $A[k]$ and the other in $B[k]$. This is usually called the
cut defined by $\vec{f}_k$. The edges $\vec{e}_j$  with this property have been directed
as going from $A[k]$ to $B[k]$.

Now consider the 1-cycle $Z[j]$ and its evaluation $<Z[j], 1 \cdot \vec{f}_k>$ ;
it equals zero unless $\vec{e}_j $ has one end in $A[k]$ and the other in $B[k]$.
It has sign $-1$ as $\vec{e}_j $ goes from $A[k]$ to $B[j]$.

Hence,
$$
 \Sigma_j \ < z[j], 1 \cdot \vec{f}_k> \ \vec{e}_j = - \Sigma_{summed \ as \ below} \  Y[j,k] \ \vec{e}_j
$$
with  the above sum is only  over the directed  edges $\vec{e}_j$ of $G-T_0$ going from $A[k]$ to $B[k]$
with direction   from $A[k]$ to $B[k]$. Note $Y[j,k] =0$ otherwise.

In total  then
$$
\Pi(\delta( \Sigma_{b \in B[k]} \ 1 \cdot b ) = 1 \cdot \vec{f}_k - \Sigma_j  \ Y[j,k] \ \vec{e}_j.
$$
Thus there is a well defined  mapping
$$ \begin{array}{l}
F: C_d(T_0;Z) \rightarrow \Pi(B^d(X;\mathbb{Z})),  \ \ \ \
\vec{f}_k \mapsto 1 \cdot \vec{f}_k - \Sigma_j \ Y[j,k] \ \vec{e}_j.
\end{array}
$$

Now let $a $ be  any element of $\Pi(B^1(X;Z))$ and note that
$$
<a, Z[j]> =0 \ as \ Z[j] \in Z_1(X;\mathbb{Z}) \ and \ \Pi(B^1(X;\mathbb{R}))^\perp = Z_1(G;\mathbb{R}).
$$
Say
$
a= ( \ \Sigma_j \ x_j \ \vec{e}_j \ )+ ( \ \Sigma_k \ y_k \ \vec{f}_k \ ).
$
Then the difference
$
Diff = a -( \ \Sigma_k  \ y_k \ F(\vec{f}_k \ ) \ )
$
lies in $\Pi(B^d(X;\mathbb{Z}))$ and is supported in $G-T_0$. But then
$
0 = <Z[e_j], Diff> = x_j \ for \ each \ j =1, \cdots, |E(G-T_0)|
$
proving that $Diff=0$. Hence, the integral chains  $F(\vec{f}_k)$ generate $\Pi(B^d(X;\mathbb{Z}))$
They are independent since projecting to $C_1(T_0;\mathbb{Z})$ each goes to $1 \cdot \vec{f}_k$,
part of the  basis
of $C_1(T_0;Z)$.

Hence, the lemma is proved.

\section{Flux and Eigenvalue Inequalities}  \label{sect8}

For the  Kirchoff Laplacian $\Delta(G)$ of a graph $G$ there is the  foundational Rayleigh-Ritz
inequality  providing an upper estimate of the first
positive eigenvalue.\footnote{An inequality giving a lower estimate on the eigenvalues is
due to Cheeger for Riemannian manifolds \cite{Cheeger1}. Bounds either way in the setting of graph theory are often referred to as Cheeger type inequalities inspired by his methods, i.e., see \cite{Chung1,Chung2}. }  Following the principle that the mesh Laplacian
is a natural generalization of the Kirchoff Laplacian, it is natural
to discuss an extension of  such inequalities to the mesh context.   The flux of 
the of the graph $G$  is a measured by the eigenvalues
of the Kirchoff Laplacian.

The mesh Laplacian $\Delta(G,T_0)= D \cdot D^\star$ is represented by $Y \cdot Y^t$
which is real symmetric with non-negative eigenvalues. Let $\Lambda$ be the
smallest positive eigenvalue of $Y \cdot Y^t$, equivalently of $\Delta(G,T_0)$. This $\Lambda$  is also the smallest positive eigenvalue  of $Y^t \cdot Y $ which represents  the reduced mesh matrix $Mesh^\#(G,T_0)$. 

By the standard Rayleigh-Ritz theorem, there is the associated equality giving upper
estimates for $\Lambda$. Namely,
$$
\Lambda = min_{f \neq 0, \ \  f \in [kernel \ C_1(T_0; \mathbb{R}) \stackrel{D^\star}{\rightarrow}
C_1(G-T_0; \mathbb{R})]^\perp  } \ \ \frac{ < D^\star(f), D^\star(f)>}{< f,f>}
$$
where $[kernel \ C_1(T_0; \mathbb{R}) \stackrel{D^\star}{\rightarrow}
C_1(G-T_0; \mathbb{R})]^\perp$ is the orthogonal complement of this kernel  in $C_1(T_0; \mathbb{R})$
endowed with the standard cellular inner product. Here $< D^\star(f), D^\star(f)> \newline = < f, D \cdot D^\star (f)>
= < f, \Delta(G,T_0)(f)>$.
Hence, the kernel consists precisely of the eigenvectors of $\Delta(G,t_0)$ with eigenvalue zero.

One may regard $\Lambda$ as a measure of the \textbf{flux} of $\Delta(G,T_0)$ through
$G-T_0$ from $T_0$ to $T_0$.

\vspace{.1in} Let the edges in $E(G-T_0)$ be enumerated and directed as before.
 Then there are three different types of directed edges.  Firstly $\vec{e}_j[P_j, Q_j]$ may be a loop, i.e.,
$P_j=Q_j$. Let these be called type 1 edges. Secondly, the edge $\vec{e}_j$ may have distinct end points,
but these end points $P_j,Q_j$ [ which lie in $V(T_0)= V(G)$] are adjacent in the subgraph $T_0$.
Call these type 2 edges. Note that in the first case $D(1 \cdot \vec{e}_j[P_j,Q_j]) = 0$ and in the second
case $D(1 \cdot \vec{e}_j[P_j,Q_j])$ is the unique directed edge in $T_0$ from $Q_j$ to $P_j$. Finally and
thirdly, the directed edge $\vec{e}_j[P_j,Q_j]$ is called of type 3, if it  has a directed simple path in $T_0$ from $Q_j$ to $P_j$
of length $\ge 2$.

In the third type 3 case, $\vec{e}_j[P_j,Q_j]$,   $P_j, Q_j$ are vertices of the
spanning tree $T_0$ so by $P_j \neq Q_j$ there is a unique directed simple
path from $Q_j$ to $P_j$ in the spanning tree $T_0$, given by distinct vertices  say $Q_j= R_{j,1}, R_{j,2}, \cdots, R_{j,n[j]}= P_j$ with $R_{j,t} \sim R_{j,t+1}$ [adjacent]  and the edge
$  [ R_{j,t} \rightarrow R_{j,t+1}]$ chosen in $T_0$. Here $n[j] \ge 2$, so one has two distinct
well defined directed edges, initial and final directed inward into $T_0$ given by
$$
[Q_j \rightarrow  R_{j,2}] = [R_{j,1} \rightarrow  R_{j,2}] \in E(T_0) \ \ and \ \ [R_{j,n[j]}\rightarrow  R_{j,n[j]-1}] =  [P_j \rightarrow  R_{j,n[j]-1}] \in E(T_0)
$$
For a directed edge $\vec{e}_j[P_j,Q_j] \in E(G-T_0)$ of type 3, let these directed two directed edges
be denoted by $ F_1( \vec{e}_j[P_j,Q_j]) =[Q_j \rightarrow  R_{j,2}] \in E(T_0)$, the initial edge of $D(\vec{e}_j[P_j,Q_j])$, and $ F_2( \vec{e}_j[P_j,Q_j]) =[P_j \rightarrow  R_{j,n[j]-1}] \in E(T_0)$, the final edge of $D(\vec{e}_j[P_j,Q_j])$ both directed inward to $T_0$.

Let $\vec{e}_j[P_j,Q_j]$ range over all type 3 directed edges of $E(G-T_0)$, set
$X$ to be the set of all  these inward directed edges so obtained:
$$
X = \{   F_1( \vec{e}_j[P_j,Q_j]) \} \cup \{ F_2( \vec{e}_j[P_j,Q_j]) \}.
$$
Naturally a given such directed edge may arise in several ways; multiplicities are neglected
in this representation. The construction provides  a set of distinct type 3 inward directed edges
of $G$.

Now form a  graph $W$ with vertices the elements of the set $X$ of  the type 3 edges.
Declare that two distinct vertices, say  with representatives $ [F_r( \vec{e}_s[P_s,Q_s])],  [F_t( \vec{e}_u[P_u,Q_u])]$,
having an edge in the graph $W$ between them  for each  type 3 directed  edge $\vec{e}_j[P_j,Q_j]$ with ingoing edges equal to these two,
i.e.,   $ \{F_1( \vec{e}_j[P_j,Q_j]),  F_2( \vec{e}_j[P_j,Q_j])\} \newline =
\{ F_r( \vec{e}_s[P_s,Q_s]),  F_t( \vec{e}_u[P_u,Q_u])\}$.

\begin{thm} \label{thmeigen} Let $\lambda$ be the smallest positive eigenvalue of the
graph Laplacian $\Delta(W)$ of the  graph $W$ derived from $G$. Let $\Lambda$ be the smallest
positive eigenvalue of the mesh Laplacian $\Delta(G,T-0)$.
Then the following inequality holds:
$$
\Lambda \le \lambda
$$
\end{thm}

By the standard methods, \cite{Chung1, Chung2} $\lambda$ has upper estimates of the Cheeger type.
These then imply upper estimates for $\Lambda$.

In detail,
let the graph $W$ have M distinct connected components, say
$$
W= \sqcup_{k=1}^M \ W[k]\ and\ correspondingly \ vertex \ decomposition \ X = V(W)
= \sqcup_{k=1}^M \ X[k].
$$
For each $k$ make a choice of proper subset
$C[k] \subset X[k]$ and denote the  complement in $X[k]$ by  $D[k] = X[k] - C[k]$,
then
$$  \begin{array}{l}
\Lambda/2 \le \lambda/2 \\
\le min \{  \frac{    \# \vec{e}_j[P_j,Q_j] \ of \ type \ 3 \ with \
                      F_1( \vec{e}_j[P_j,Q_j]) \in C[k] \ and \  F_2( \vec{e}_j[P_j,Q_j]) \in D[k]
                     \  or \ visa \ versa   }{  min(|A[k]|,|B[k])  } , k=1, \cdots, M  \}
                       \end{array}
                     $$

   Proof: Let $j: R[X] \subset C_1(T_0; \mathbb{R})$ denote the inclusion of the real sums of
   multiples of the [inward directed]  edges in $X$. Let $\overline{D^\star}= D^\star \cdot j : R[X] \rightarrow C_1(G-T_0; \mathbb{R})$
   with adjoint $(\overline{D^\star})^\star$. Let $\mu$ be the smallest positive eigenvalue
   of the composite $\overline{\Delta} = (\overline{D^\star})^\star \cdot \overline{D^\star}$. Tautologously,
   one has the inequality [restricting the minimum over fewer cases]:
   $$
   \Lambda \le  min_{g \neq 0, \ \  g \in R[X] \cap [kernel \ C_1(T_0; \mathbb{R}) \stackrel{D^\star}{\rightarrow}
C_1(G-T_0; \mathbb{R})]^\perp  } \ \ \frac{ < D^\star(g), D^\star(g)>}{< g,g>}= \mu
$$

Now here $R[X] \cap [kernel \ C_1(T_0; \mathbb{R}) \stackrel{D^\star}{\rightarrow}
C_1(G-T_0; \mathbb{R})] $  consists of the sums $Sum= \Sigma_{x\in X} \ A[x] \ x$
for which $D^\star(Sum)=0$. That is, for each $\vec{e}_j[P_j,Q_j]$ the equation \newline
$< D(\vec{e}_j[P_j,Q_j]),Sum> = 0$ holds. But $< D(\vec{e}_j[P_j,Q_j]),Sum> =  < F_2( \vec{e}_j[P_j,Q_j]) - F_1(\vec{e}_j[P_j,Q_j]), Sum> = A[F_2( \vec{e}_j[P_j,Q_j]) ] -A[F_1(\vec{e}_j[P_j,Q_j])]$. This means that if the vertices in $W$,
$x,y$ are adjacent, then $A[x] = A[y]$.  Consequently, $R[X] \cap [kernel \ C_1(T_0; \mathbb{R}) \stackrel{D^\star}{\rightarrow}
C_1(G-T_0; \mathbb{R})]$ consists of the sum of multiples of the entries one
for each component
$$
S[k] =( \  \Sigma_{x \in W[k]}  \ 1 \cdot x \ ) \in C_1(T_0; \mathbb{R}), \ k=1, \cdots, M
$$
and the mapping $(\hat{D^\star})^\star: C_1(G-T_0; \mathbb{R}) \rightarrow \{ \Sigma_{x \in X} \ A[x] \ x\}$
sends $ 1 \cdot \vec{e}_j(P_j,Q_j) =  F_1(\vec{e}_j[P_j,Q_j]) - F_2(\vec{e}_j[P_j,Q_j])$.

That is, under this identification, $(\overline{D}^\star)^\star$ is identified with
$\partial :C_1(W; \mathbb{R}) \rightarrow C_0(W; \mathbb{R})$ and $\overline{D}^\star$ with $\partial^\star: C_0(W; \mathbb{R}) \rightarrow C_1(W; \mathbb{R})$ and the kernel of $\partial^\star$ with the sums of multiples of the $S[k]$'s.
The orthogonal complement of this kernel is identified  with the span of the eigenspaces with positive eigenvalues.
In particular, $\mu = \lambda$ is given as the least positive eigenvalue
of the graph Laplacian of the graph $W$ in the above inequality.

The standard Cheeger inequality \cite{Chung1,Chung2} now applies to the graph $W$ to give the above estimate.

\section{Appendix; A Correction to the statement of theorem 7.1 of \cite{Cappell1}
and extension of theorem \ref{thm1} to the higher dimensional context :} \label{sect9}

The paper \cite{Cappell1} considers  generalizations of the mesh matrix
to the context of $d$-dimensional finite CW-complexes $X$. The theorem 7.1
of \cite{Cappell1}, once stated correctly,  is a generalization of lemma \ref{lemmaA}.
Once lemma \ref{lemmaA} is generalized, the methods of this paper combined with
\cite{Cappell1} also give a generalization
of theorem \ref{thm1} above. This  is carried out in this section.
In the context of the higher dimensional Kirchoff Laplacian for CW-complexes
some relevant papers are \cite{Adin1,Bernardi1,Bolker1,Farid1, Klein1, Klein2,Klein3, Duval1, Duval2, Kalai1, Lyons1,Lyons2,Peterson1} and others.

Let $X$ be a finite d-dimensional CW-complex. Let $S_d$ be the
$d$-cells of $X$ with chosen orientations, i.e., $X$ is obtained from the $d-1$ skeleton
of $X$, i.e., $X^{d-1}$,  by attaching the oriented $d$-cells of $S_d$.   For $V \subset S_d$ define  the subcomplex
$$
X_V = X^{d-1} \cup V \subset X^{d-1} \cup S_d = X.
$$
For a finite cell complex, say $X$, let $t_d(X)$ denote the order of the torsion
subgroup of the $(d-1)^{th}$ integral homology of $X$, $H_{d-1}(X;\mathbb{Z})$.

In this higher dimensional context a  subset $V \subset S_d$
is called a \textbf{spanning forest} of $X$ if the composition
$$
C_d(X_V;\mathbb{R}) \stackrel{\partial_d}{\rightarrow} B_{d-1}(X_V;\mathbb{R}) \subset B_{d-1}(X;\mathbb{R})
$$
is an isomorphism. In the case that  $G$ is a graph and  the graph $G$ is connected,
the spanning forests in the above sense are exactly the spanning trees of the graph $G$.
In the case that a graph $G$ has $K$ connected components, say $G = \sqcup_{j=1}^K \ G[j]$,
a spanning forest is a choice for each connected component $G[j]$ of  a spanning tree $T[j]$ of $G[j]$.
It is shown in \cite{Cappell1} that all $d$-dimensional CW-complexes have spanning forests.

Now the integral $d$-cycles $Z_d(X;\mathbb{Z})$ is a free abelian group, so one
may pick a integral basis, say $\{z_j\}$. Via the natural pairing
$$
Z_d(X;\mathbb{Z}) \times Z_{d}(X;\mathbb{Z}) \subset C_d(X;\mathbb{Z}) \times C_d(X;\mathbb{Z})
\stackrel{<.,.>}{\rightarrow} \mathbb{Z},
$$
one may form the mesh matrix for this integral choice of basis
$$
Mesh(X; \{z_r\}) = \{ < z_a, z_b> \}
$$

\vspace{.3in}
On the other hand, if one chooses a spanning forest, say $V_0$, of $X$, for each
$e$ in $ S_d - V_0$, the boundary $\partial_d (e)$ is \textbf{uniquely expressed} as $\partial_d ( D(e))$
for a unique $d$-chain $ D(e) \in C_d(X_{V_0}; \mathbb{R})$. In particular, for each $e\in S_d-V_0$
there the unique $d$-cycle
$$
z(e) = 1 \cdot e - D(e) \in Z_d(X; \mathbb{R})
$$
with $D(e)$ supported exactly in $X_{V_0}$. Since $Z_d(X_{V_0}; \mathbb{R}) =0$ by $V_0$
a spanning tree, one sees that these $\{z(e) \ | \ \ e \in S_d - V_0\}$
form a ``geometric'' basis for $Z_d(X;\mathbb{R})$ with associated mesh matrix
$$
Mesh(X;\{ z(e) \ | \ e \in S_d - V_0\} , geometric)
$$
Note that $z(e)$ is supported in $X_{e}= X^{d-1} \cup e$ for $ e \in S_d-V_0$.
Let $N= |S_d-V_0|$.

Let the $d$-cells in $S_d - V_0$, be enumerated via $e_1, \cdots, e_{N}$.
Consider  a choice of $j \times j$ diagonal minor, say $Mesh(\{z(e)\}, geometric)_{k_1,k_2, \cdots, k_j}$
with $ 1 \le k_1 <k_2< \cdots <  k_j \le N$ by picking these $j$ rows and columns. The
entries of this minor all are of the form
$< z(e_{k_a}), z(e_{k_b})> $ with each
 $z(e_{k_a}) \in
Z_d(X_{V_0 \cup \{e_{k_1},e_{k_2}, \cdots, e_{k_j}\}}  ); \mathbb{R})$ for $1 \le 1 \le j$.
Consequently,  as in \S 3,
$$
det(Mesh(\{X(e)\}), geometric)_{k_1,k_, \cdots, k_j})= det( Mesh(X_{V_0 \cup \{e_{k_1}, \cdots, e_{k_j}\}}), geometric).$$
Consequently, modeled on  the proof of lemma \ref{lemmaA}, one has for the geometric basis the equality
$$
det( X \cdot Id - Mesh(X; \{z(e)\ | \ in \ S_d-V_0\}, geometric)
= X^N + \Sigma_{j=1}^N \ (-1)^j  \ b_j \ X^{N-j}
$$
with $b_j = \Sigma_{1 \le k_1<k_2, \cdots , k_j \le N} \ \  det( Mesh(  X_{V_0 \cup \{e_{k_1}, \cdots, {e_{k_j}\}}})),
geometric))$. Here the choice of $k_1,k_2, \cdots, k_j$ in $S_d-V_0$ may be
recorded as a choice of subset $U=\{k_1, \cdots, k_j\} \subset S_d - V_0$ with $|U|=j$.

 \vspace{.3in}
 In these terms    the slightly corrected  version of theorem 7.1 of \cite{Cappell1} is:
[It is the high dimensional version of lemma \ref{lemmaA}.]
\begin{thm}  \label{thmcorrectedversion}
Let $V_0 \subset S_d$ be a $d$-dimensional spanning forest of a $d$-dimensional CW-complex $X$
with oriented $d$-cells $S_d$, then $\sigma_j(X;V_0) = $
the coefficient of $(-1)^j \ T^{N-j}$ in the characteristic polynomial
$$
det( Mesh(X, \{ z(e) \ | \ e \in (S_d-V_0) \}), geometric  \ basis)
$$
is given by
$$
\sigma_j(X;V_0) = \Sigma_{U \subset (S_d-V_0) \ with \ |U| = j}
\ \ ( \  \Sigma_{V \ a \ spanning \ forest \ of\  X_{V_0 \cup U} } \ \left( \frac{t_{d-1}(X_V)  }{t_{d-1}(V_0) } \right)^2 \ )
 $$
\end{thm}

For $d=1$, there is no torsion in the homologies, so sum is over entries all equal to $+1$.
So, this  is the high dimensional version of lemma \ref{lemmaA}. The proof is a consequence of  the above argument
 and the following equality: For $Y$  a d-dimensional finite CW-complex, and a choice of spanning forest $V_0 \subset S_d$ for $Y$.
 $$
 det( Mesh(Y, \{z(e)\}, geometric) =
   \Sigma_{V \ a \ spanning \ forest \ of\  Y } \ \left( \frac{t_{d-1}(Y_V)  }{t_{d-1}(Y_{V_0}) }  \right)^2. \ \ (\star)
 $$
 applied to the various $Y= X_{V_0 \cup U}$.

 This last is the analog of Trent's theorem \ref{thmTrent1} above. In the paper \cite{Cappell1}
 theorem 2.1 asserts the evaluation of the determinant of the mesh matrix for the
 \textbf{integral basis} $\{z_r\}$ as:
 $$
 det( Mesh(X;\{z_r\})) =   \Sigma_{V \ a \ spanning \ forest \ of\  X } \ \left( \frac{t_{d-1}(X_V)  }{t_{d-1}(X) } \right)^2   \ \ ( \star \ \star)
 $$

 This equality $( \star \ \star)$ is to be modified to get $(\star)$ which utilizes the geometric basis
 rather than the geometric basis.

 Here under the mapping $Z_d(X;\mathbb{R}) \subset C_d( X; \mathbb{R})
  \rightarrow C_d(X, X_{V_0}; \mathbb{R})$  [which by $V_0$ a spanning forest
  is an isomorphism] the lattice of the geometric basis $Z\{z(e) \ | \ e \in S_d-V_0\}$
   has image equal to the image of the integral classes $C_d(X;\mathbb{Z})$, [as
    $z(e) = 1 \cdot e + D(e)$ with $D[e] \in C_d(X_{V_0}; \mathbb{R})$]; and the lattice of the
  integral  basis for $Z_d(X;\mathbb{Z})$ has image equal of course  to the image of the integral
 basis of  $Z_d(X;\mathbb{Z})$. As $V_0$ is a spanning forest, these images provide
  two real bases for the $C_d(X, X_{V_0}; \mathbb{R})$ with
   one lattice inside the other. Hence, the change of basis
  from integral to geometric
  has determinant equal to the order of the finite abelian group:
  $$
   \frac{C_d(X, X_{V_0}; \mathbb{Z})}{ image \ Z_d(X;\mathbb{Z})}.
  $$
  The order of this finite abelian group is identified on page 14 of \cite{Cappell1} as
  $$
  | \frac{  C_d(X, X_{V_0}; \mathbb{Z})   }{   Z_d(X;\mathbb{Z})}| = \frac{t_{d-1}(X_{V_0})}{t_{d-1}(X)}
  $$
  The respective determinants of mesh matrices are related   by the square of this quotient.
  This relation demonstrates that $(\star \ \star)$ implies $(\star)$ so the
  above forms  an independent proof of theorem \ref{thmcorrectedversion}.

 \vspace{.3in}

  Now as above in the present paper, the geometric basis is of the form:
  $z(e_j) = 1 \cdot e_j + D( e_j)\in Z_d(X;\mathbb{R})$ for $e_j \in S_d-V_0$ with $D( e_j)$ supported
in $X_{V_0}$. This implies
$$
<z(e_j), z(e_k)> = Id_{j,k}+ < D(e_j), D(e_k)>.
$$
That is, the geometric mesh matrix is of the form,
$$
Mesh(X;\{ z(e_j)\}, geometric) = Id + Mesh^\#(X; \{z(e_j)\}, geometric)
$$
with the reduced mesh matrix $Mesh^\#(X; \{z(e_j)\}, geometric)$  defined by \newline $\{< D(e_j), D(e_k)>\}$.
Hence, the counting argument given for proving theorem \ref{thm1} of this paper applies
without essential change to  yield
 a higher dimensional generalization of this paper's theorem \ref{thm1}. Namely:

\begin{thm}  \label{thmhigher}
$$ \begin{array}{l}
det( (U+1) \ Id - Mesh(X; \{z(e_a)\}, geometric)) \\ = det( U \ Id - Mesh^\#(X; \{z(e_a)\}, geometric)
 = U^{N} + \Sigma_{j=1}^{N} \ (-1)^j \  c_j  \  U^{N-j} \\
with \\
c_j =  \Sigma_{U \ a \ spanning \ forest \ of\  X\ with \ |U \cap (S_d-V_0)| = j  } \ \left( \frac{t_{d-1}(X_U)  }{t_{d-1}(X_{V_0}) } \right)^2
\end{array}
$$
\end{thm}

\vspace{-0.1in}

\frenchspacing

\bibliographystyle{plain}

Sylvain E. Cappell,  cappell@courant.nyu.edu, Courant Institute of Mathematical Sciences, NYU.

Edward Y. Miller, emiller@cims.nyu.edu, , Courant Institute of Mathematical Sciences, NYU.

\end{document}